\documentclass[letterpaper, 10 pt]{article}

\usepackage{graphicx}
\usepackage[dvipsnames]{xcolor}
\usepackage{soul}
\usepackage{url}
\usepackage[algo2e]{algorithm2e}
\RestyleAlgo{ruled}
\usepackage{subcaption}
\makeatletter
\let\NAT@parse\undefined
\makeatother
\usepackage{hyperref}
\usepackage{tikz}
\usepackage{bbm}
\usepackage{epsfig}
 
\usepackage{amssymb}
\usepackage{amsthm}
\usepackage{amsmath}
\usepackage{cancel, comment, nicefrac}
\usepackage{mathtools}
\usepackage[nolist]{acronym}

\usepackage[shortlabels]{enumitem}

\newtheorem{theorem}{Theorem}
\newtheorem{proposition}{Proposition}
\newtheorem{lemma}{Lemma}
\theoremstyle{definition}
\newtheorem{definition}{Definition}
\newtheorem{assumption}{Assumption}
\newtheorem{remark}{Remark}

\definecolor{matt}{rgb}{0.5 0 1}
\setstcolor{matt}

\newcommand{\argmin}[1]{\underset{#1}{\mathrm{argmin}}}

\newcommand{\tr}{\mathrm{tr}}

\newcommand{\cD}{\mathcal{D}}

\newcommand{\cG}{\mathcal{G}}

\newcommand{\cM}{\mathcal{M}}

\newcommand{\cQ}{\mathcal{Q}}

\newcommand{\cU}{\mathcal{U}}

\newcommand{\cX}{\mathcal{X}}

\newcommand{\E}{\mathbb{E}}

\newcommand{\N}{\mathbb{N}}

\newcommand{\R}{\mathbb{R}}

\usepackage[margin=2.5cm]{geometry} 

\begin{acronym}
    \acro{OCP}{Optimal Control Problem}
    \acro{HJB}{Hamilton-Jacobi-Bellman}
    \acro{LP}{Linear Program}
    \acro{LMI}{Linear Matrix Inequality}
    \acro{CD}{Christoffel-Darboux}
    \acro{MPC}{Model Predictive Control}
    \acro{SoS}{Sum-of-Squares}
\end{acronym}

\title{\bf Convex computation of cost bounds for data-driven nonlinear optimal control}
\author{Oumayma Khattabi$^{1}$, Matteo Tacchi-Bénard$^{2}$, Sorin Olaru$^{1}$
\thanks{The authors would like to thank Colin Jones for fruitful discussions.}
\thanks{This research was supported by the French government under the France 2030 program, reference ANR-11-IDEX-0003 within the OI H-Code, by doctoral school STIC, and by graduate school ISN of Paris-Saclay University.
\newline The work of Oumayma Khattabi and Sorin Olaru is supported by the RTE Chair at CentraleSupelec.
\newline The work of Matteo Tacchi-Bénard is supported by the French National Research Agency (ANR) under the project grant ANR-25-CE48-1599-01.
\newline $^{1}$ Université Paris-Saclay, CNRS, CentraleSupélec, Laboratoire des Signaux et Systèmes (L2S), France.
\newline $^{2}$ Univ. Grenoble Alpes, CNRS, Grenoble-INP (Institute of Engineering Univ. Grenoble Alpes), GIPSA-lab, 38000 Grenoble, France.
}}

\begin{document}

\maketitle

\begin{abstract}
The paper addresses the optimal control problem of an input-affine polynomial system with an unknown (state-dependent) input matrix affecting the control effectiveness. Under assumptions on the existence of a data set and the Lipschitz continuity of the unknown part  {of the dynamics}, an uncertainty set is constructed and used to adapt the Hamilton-Jacobi-Bellman equation towards a data-driven formulation. We solve it for the best-case scenario, thus giving a lower bound on the best cost achievable by means of a feedback control policy. After relaxing it into a moment-SoS hierarchy of semidefinite programs, a control law is extracted from the primal on moments and compared with other benchmark methods in a numerical application.

\end{abstract}

\section{Introduction}

We consider the problem of optimal control design for input-affine systems. The control-effectiveness matrix field is treated as unknown or partially known, and uncertainty is reduced through state-space evaluations and bounded by Lipschitz continuity properties. Such a problem draws from several interconnected domains: nonlinear optimal control and \ac{HJB} formulations, set-membership/bounded-error estimation with Lipschitz bounds, data-driven control, and learning-based and myopic control of unknown dynamics.

Optimal control for nonlinear systems is a well-trodden problem in control theory~\cite{dontchev1995optimality}, and many methods have been developed to solve it when the model is available. One can cite direct shooting methods~\cite{bock1984multiple}, differential dynamic programming~\cite{mayne1966second}, and nonlinear \ac{MPC}~\cite{mayne2000constrained,diehl2005real} among others. There are also approaches solving the \ac{HJB} equation through convex reformulations as an alternative to numerical viscosity solutions~\cite{fleming2006controlled,bardi1997optimal}. Among these, one finds~\cite{vinter1993convex}, where the \ac{OCP} is formulated as an infinite-dimensional \ac{LP} over occupation measures, and its dual yields subsolutions to the \ac{HJB} equation. Later, \cite{lasserre2008nonlinear,henrion2008nonlinear} apply this formulation to polynomial systems.
Global feedback laws are extracted in the same setting~\cite{majumdar2014convex}, with convergence guarantees in~\cite{korda2016controller} and performance bounds in~\cite{jones2020polynomial}.

When part or all of the model is unavailable, as is the case with many systems, data-driven  {control} approaches become necessary. Some methods recover a model from data~\cite{BPK16}, and other works use such learned models to generate control laws via \ac{MPC}, for example, as in~\cite{KKB18,MLMC20}. In other approaches, the uncertainty is bounded in a set membership  {framework}, where all models consistent with the data are identified~\cite{milanese2004set}, or used by a stable estimator to bound state and actuator faults~\cite {buciakowski2017bounded}. These bounds are also used to provide the control law directly, as is the case with myopic control in~\cite{ornik2019myopic} and data-driven set-membership \ac{MPC}~\cite{wei2026data}. Finally, there are approaches in which the model and bounds are bypassed, so control is performed directly from data, as in classic data-enabled predictive control~\cite{coulson2019data}.

In this work, we adapt the model-based \ac{OCP} described in~\cite{lasserre2008nonlinear}, restricted to control-affine systems, to handle a partially unknown model. The particularity of the optimal control design in the presence of uncertainty in the affine term lies in the fact that the equivalent differential inclusion loses its convexity properties, among other strong properties of the exact model-based counterparts. Consequently, the pairs of states and inputs have to be analyzed in terms of both distance in the state space and excitability with respect to input values. The corresponding bounds on the optimal data-driven cost need to be evaluated. The main contributions of this work are a ``best-case'' lower bound on the attainable cost, a guarantee of convergence for the data-driven \ac{OCP}, a data selection scheme, and an alternative polynomial control law derived from occupation measures, whose performance is then compared to \ac{MPC}-based approaches~\cite{KKB18,MLMC20}.

The present paper is organized as follows. Section~\ref{sec:problem_setting} provides the problem setting and assumptions on the system. Section~\ref{sec:convergence} presents the data-driven \ac{OCP} and the convergence proof. In Section~\ref{sec:msos}, the moment-\ac{SoS} hierarchy framework is detailed. Section~\ref{sec:extract} shows how the control law is extracted. Numerical applications where performance is compared to that of benchmark methods are detailed in Section~\ref{sec:num_app}.

\vspace{.2cm}
\noindent{\bf Notation:} Let $n,d \in \N$, $[n] \triangleq \{1,\ldots,n\}$.
$\|.\|$ is the Euclidean norm for vectors, spectral norm for matrices, and the absolute value for scalars.
$\|.\|_F$ is the Frobenius norm for matrices. 
$\E[x]$ denotes the expected value of the random variable $x$, $\tr(A)$ the trace of matrix $A$. 
$\R[x]^n$ is the set of $n$-dimensional vectors of real polynomials in $x\in\R^n$, and $\R_d[x]$ is the subset with polynomials up to degree $d$.
$C^n(X)$ is the set of functions on $X$ with continuous derivatives up to order $n$. $\cM(X)$ is the set of nonnegative Borel measures over $X$.

\section{Data-driven optimal control}\label{sec:problem_setting}

Consider the nonlinear \ac{OCP}
\begin{subequations}\label{eq:OCP}
\begin{align}
    V^\star \triangleq \inf_{u(\cdot):[0,T] \to \cU}& \int_0^T \ell(x(t),u(t)) \; dt + h(x(T)) \\
    \text{s.t.} \quad \; & \quad \dot{x}(t) = f(x(t)) + G(x(t)) \, u(t) \label{eq:sys} \\
    & \quad x(t) \in \cX, \quad x(T) \in \cX_T \label{eq:set_constraint}
\end{align}
\end{subequations}
where  { $T>0$ is the time horizon,} $f,h,\ell$ are polynomials and $G: \cX \longrightarrow \R^{n\times m}$ is \textit{unknown}.  {The sets  $\cX$ and $\cX_T$ in the state space} are semi-algebraic and compact:
$$ \cX \triangleq \{x \in \R^n \mid p_\cX(x) \geq 0\}, \ \ \cX_T \triangleq \{x \in \R^n \mid p_T(x) \geq 0\} $$
with $p_\cX, p_T$ vectors of polynomials (inequalities are meant component-wise), and the admissible input set $\cU \subset \R^m$ is described by means of a convex polytope:
$$ \cU \triangleq \{u \in \R^m \mid p_\cU(u) \triangleq P_\cU \, u - q_\cU \geq 0 \}, $$
with $P_\cU \in \R^{D \times m}$ and $q_\cU \in \R^D$. The objective of this work is to construct a control input by solving an \ac{OCP}, accounting for the uncertainty on $G(\cdot)$, under the following assumptions.
\begin{assumption}\label{ass:dataset}
A dataset of finite and noise-free evaluations $\cD\triangleq\{(x_i,u_i,\dot{x}_i)\}_{i\in[N]}\subset \cX \times (\cU\setminus\{0\})\times \R^n$, is available.
\end{assumption}

\begin{assumption}\label{ass:Lipschitz}
The unknown matrix function $G(\cdot)$ is Lipschitz continuous on the set $\cX$, and an upper bound $K>0$ on its Lipschitz constant is known:
$$\|G(x) - G(x')\| \leq K \|x-x'\| \qquad \forall x,x' \in\cX.$$
\end{assumption}

For $i\in[N]$, let $y_i \triangleq \dot{x}_i-f(x_i) = G(x_i) \, u_i$. Based on assumptions~\ref{ass:dataset} and~\ref{ass:Lipschitz}, we get
\begin{equation}\label{eq:Lipschitz}
    \|G(x) \, u_i - y_i\| \leq K \|u_i\| \|x-x_i\| \qquad \forall x \in\cX.
\end{equation}

By normalizing inequality~\eqref{eq:Lipschitz} over  {$\|u_i\|\neq 0$}, we define the set containing all possible values of $G(x)$ as the set-valued map
\begin{equation}\label{eq:uncertainty_set}
    \cG(x) \triangleq \left\{ B \in \R^{\scriptscriptstyle n \times m} \ \mid \ p_i(x,B) \geq 0 , \ \forall i \in [N] \right\},
\end{equation}
with $p_i(x,B) \triangleq K^{\scriptscriptstyle 2} \|x - x_i\|^{\scriptscriptstyle 2} - \left\| B\frac{u_i}{\|u_i\|} - \frac{y_i}{\|u_i\|} \right\|^{\scriptscriptstyle 2}$. Then, we deduce the global \emph{semi-algebraic} uncertainty set
\begin{equation}\label{eq:Gamma}
    \!\!\!\! \Gamma \triangleq \left\{ (x,B) \in \cX \times \R^{\scriptscriptstyle n  \times  m} \ \mid \ p_i(x,B) \geq 0 , \ \forall i \in [N] \right\}.\!
\end{equation}

\begin{remark}
    The semi-algebraic uncertainty envelope $\Gamma$ is bounded as soon as the input data $u_i$ span $\R^m$.
\end{remark}

\begin{remark}
The definition of the uncertainty set $\cG$ does not exclude non-Lipschitz behavior from being included in $\Gamma$, as long as the Lipschitz inequality is verified on the data points. The set $\Gamma$ is less restrictive than Assumption~\ref{ass:Lipschitz}, yet such a relaxed account of the assumptions is still permissible in the proposed approach.
\end{remark}

\begin{remark}\label{rem:noise}
Incorporating bounded measurement noise into $\Gamma$ is feasible.
For $i\in[N]$, taking $\Tilde{x}_i,\dot{\Tilde{x}}_i,\Tilde{y}_i$ as the measured values, let $\varepsilon_x,\varepsilon_{\dot{x}}>0$ be such that $\|x_i-\Tilde{x}_i\| \leq \varepsilon_x$ and $\|\dot{x}_i-\dot{\Tilde{x}}_i\| \leq \varepsilon_{\dot{x}}$.
Since $f(\cdot)$ is polynomial, it is locally Lipschitz in $\cX$. Let $K_f$ be its local Lipschitz constant. For $x\in\cX$:
\begin{align*}
    \|G(x) \, u_i \! - \! \Tilde{y}_i\| = \|G(x) \, u_i \! - \! y_i + \dot{x}_i \! - \! \dot{\Tilde{x}}_i - (f(x_i) \! - \! f(\Tilde{x}_i))\| &\leq K \|u_i\| \|x-x_i\| + \varepsilon_{\dot{x}} + K_f \, \varepsilon_{x} \\
    &\leq K \|u_i\| \|x-\Tilde{x}_i\| \! + \! K \|u_i\| \varepsilon_{x} \! + \! \varepsilon_{\dot{x}} \! + \! K_f \, \varepsilon_{x}
\end{align*}
By normalizing over $\|u_i\|$ and squaring both sides, we get
\begin{gather*}
    \left\| G(x) \frac{u_i}{\|u_i\|} \! - \! \frac{\Tilde{y}_i}{\|u_i\|} \right\|^{\scriptscriptstyle 2} \leq \big(K \|x \! - \! \Tilde{x}_i\| + \Delta\big)^{\scriptscriptstyle 2}\\
\end{gather*}
with $\Delta = \big(K + \frac{K_f}{\|u_i\|} \big) \varepsilon_x + \frac{\varepsilon_{\dot{x}}}{\|u_i\|}$. For any $\theta > 0$, defining $\alpha = (1+\theta^{-1})\Delta^2$ and using $$(a+b)^2\le(1+\theta)a^2+(1+\theta^{-1})b^2,$$ $p_i$ is replaced by the
polynomial
$$\Tilde{p}^\theta_i(x,B)\triangleq(1+\theta)K^2\|x-\Tilde{x}_i\|^2+\alpha
-\left\|B\frac{u_i}{\|u_i\|}-\frac{\Tilde{y}_i}{\|u_i\|}\right\|^2,$$
which still defines a set containing $G(x)$ for all $x\in\cX$.
\end{remark}

\begin{definition}\label{def:admissible}
Let $\cU_T$ be the set of vector control functions from $[0,T]$ to $\cU\in\R^m$. An admissible control $u(\cdot)\in\cU_T$ for an initial point $x(0) = x_0 \in\cX$ is a measurable control function for which the trajectory $x(\cdot)$, solution to~\eqref{eq:sys}, stays in $\cX$ and is continuous over $[0,T]$, and $x(T)\in\cX_T$ for $\cX_T$ a compact subset of $\cX$ taken as the target set.
\end{definition}

We seek to construct an admissible control law solving (at least approximately) problem~\eqref{eq:OCP}.
One approach to this \ac{OCP} is to solve the \ac{HJB} partial differential equation
\begin{equation}\label{eq:HJB}
    \frac{\partial V}{\partial t} + \min_{u\in\cU} \left(\frac{\partial V}{\partial x} (f(x) + G(x) \, u) + \ell(x,u)\right) = 0 \hspace{-.5em}
\end{equation}
with $V:[0,T]\times\cX\rightarrow\R\cup\{+\infty\}$ being the Bellman value function verifying the terminal constraint $V(T,x) = h(x)$. For a feasible problem, the minimum $V^\star$ of~\eqref{eq:OCP} is $V(0,x_0)$.

To account for the lack of information on $G(\cdot)$, it is imperative to use the uncertainty set $\cG(\cdot)$. To take into account all possible values, $G(\cdot)$ may be replaced by an arbitrary $B\in\cG(\cdot)$, and equation~\eqref{eq:HJB} should be relaxed accordingly. Two options seem natural:  robust control accounting for the worst case, and optimistic control for the best case.
While robust control seems to be the most attractive approach in terms of applications and guarantees, it introduces a classical game-theoretical minimax structure of the Hamilton-Jacobi-Isaacs type, with the uncertain variable $B$ playing against the controller $u$. This results in intricate nonconvex bilinear behaviors that are beyond the scope of this contribution.

\section{Infinite dimensional linear programs}\label{sec:convergence}

\subsection{Generalized Moment Problems}

In this paper, we are interested in the best case, where we solve for the lowest cost possible achievable by the system within the bounds of $\Gamma$. In this case, $B$ assists the control $u$, and we obtain a lower-bound benchmark against which no other control can achieve a better outcome. The \ac{HJB} equation in~\eqref{eq:HJB} is recast to search for a best-case value function $W$:
\begin{equation}\label{eq:DD-HJB-best}
    \frac{\partial W}{\partial t} + \min_{\substack{u\in\cU \\ B\in\cG(x)}} \left(\frac{\partial W}{\partial x} (f(x) + B \, u) + \ell(x,u) \right) = 0. \hspace{-.3em}
\end{equation}

Next, following the approach in~\cite{vinter1993convex,lasserre2008nonlinear}, we relax~\eqref{eq:DD-HJB-best} into inequalities, and solve for a polynomial subsolution $w \in \R[t,x]$ of~\eqref{eq:DD-HJB-best} through the following \ac{LP}
\begin{subequations}\label{eq:LP-dual}
\begin{align}
    \hspace{-2em} w_\infty^\star \triangleq &\sup_{w\in \R[t,x]} w(0,x_0) \\
    \text{s.t.} \ \ \ &\forall t\in[0,T], \ (x,B)\in\Gamma, \ u\in\cU, \notag
\end{align}
\vspace{-1.5em}
\begin{align}
    \frac{\partial w}{\partial t}(t,x) \! + \! \frac{\partial w}{\partial x}(t,x) \Big(f(x) \! + \! B \, u\Big) & \geq - \ell(x,u) \label{eq:subHJB} \\
    \forall x \in \cX_T, \qquad w(T,x) &\leq h(x) \label{eq:subTerminal}
\end{align}
\end{subequations}

\begin{proposition}\label{prop:inequalities}
Let $w_\infty^\star$ be the optimal value of~\eqref{eq:LP-dual}, and $W$  the unique solution to~\eqref{eq:DD-HJB-best} corresponding to $W(T,x) = h(x)$. Then $w_\infty^\star \leq W(0,x_0) \leq V(0,x_0) = V^\star \leq J(u(\cdot))$, for any $u(\cdot)\in\cU_T$ feasible in~\eqref{eq:OCP}, where 
$$J(u(\cdot)) = \int_0^T \ell(x(t),u(t)) \; dt + h(x(T)) \quad \text{s.t. }\eqref{eq:sys},~\eqref{eq:set_constraint}$$
\end{proposition}
\begin{proof}
First, we prove the inequality on the left-hand side. Let $\overline{x}(\cdot)$ be a trajectory solution to~\eqref{eq:sys} and~\eqref{eq:set_constraint}, corresponding to $\overline{u}(\cdot)\in\cU_T$, an admissible control for initial condition $x_0$, and $\overline{G}(\overline{x}(t))\in\cG(\overline{x}(t))$. Since $w\in\R[t,x]\subset C^1$:
\begin{gather*}
    w(T,\overline{x}(T)) - w(0,x_0) = \int_0^T \frac{dw}{dt}(t,\overline{x}(t)) dt = \int_0^T \!\! \frac{\partial w}{\partial t} + \frac{\partial w}{\partial x} \Big(f(\overline{x}(t)) + \overline{G}(t) \, \overline{u}(t)\Big) dt \overset{\eqref{eq:subHJB}}{\geq} - \int_0^T \!\!\!\! \ell(\overline{x}(t),\overline{u}(t)) dt \\
    \Longrightarrow w(0,x_0) \leq \int_0^T \!\!\!\! \ell(\overline{x}(t),\overline{u}(t)) dt + w(T,\overline{x}(T)) \overset{\eqref{eq:subTerminal}}{\leq} \int_0^T \!\!\!\! \ell(\overline{x}(t),\overline{u}(t)) dt + h(\overline{x}(T)).
\end{gather*}
With $W(0,x_0) \! = \! \inf_u \int_0^T \!\! \ell(x(t),u(t)) dt + h(x(T))$ for any feasible $(x(\cdot),u(\cdot))$, including $(\overline{x}(\cdot),\overline{u}(\cdot))$, and $w_\infty^\star = \sup_{w} w(0,x_0)$, we get $w_\infty^\star \leq W(0,x_0)$.

For the second inequality, let $V$ be the unique solution to~\eqref{eq:HJB} corresponding to $V(T,x) = h(x) = W(T,x)$. Since the real $G(x)\in\cG(x)$ for all $x\in\cX$, then $W$ is an infimum over a larger set that includes the set over which $V$ is solved; hence $W(0,x_0) \leq V(0,x_0)$. The last inequality is the definition of the optimality $V^\star = \inf J$.
\end{proof}

Problem~\eqref{eq:LP-dual} is an~\ac{LP} and as such it is subject to Lagrange duality, with the following Generalized Moment Problem on occupation measures:

\begin{subequations}\label{eq:LP-primal}
\begin{gather}
    \hspace{-2em} w_\infty' \triangleq \inf_{\substack{\mu \in \cM(\Xi) \\ \nu \in \cM(\cX_T)}} \int_{\Xi} \ell \; d\mu + \int_{\cX_T} h \; d\nu \\
    \text{s.t. } \  \delta_T\otimes\nu + \left[\frac{\partial\mu}{\partial t} + \frac{\partial}{\partial x}\Big((f+B \, u) \; \mu\Big) \right]^{(t,x)} \hspace{-1em} = \delta_{(0,x_0)} \label{eq:liouville}
\end{gather}
\end{subequations}
with $\Xi \triangleq [0,T]\times\Gamma\times\cU$. $\nu$ is called the terminal measure, while $\mu$ is called the occupation measure and encapsulates the trajectories $\xi = (t,x(t),B(t),u(t)) \in \Xi$. Constraint~\eqref{eq:liouville} is called the Liouville equation, and models how the controls transport the initial condition $x_0$ onto the distribution $\nu$ of terminal conditions.

\begin{proposition}\label{prop:strong_duality}
    Strong duality holds, i.e., $w_\infty^\star = w_\infty'$.
\end{proposition}
\begin{proof} The following choice for $w$ is strictly feasible in~\eqref{eq:LP-dual}
    $$w(t,x) = \Big(1-\min_{\substack{\tilde{x}\in\cX \\ u\in\cU}}\ell(\tilde{x},u)\Big)(t-T) + \min_{\tilde{x}\in\cX}h(\tilde{x}) - 1$$
    so Slater's condition~\cite{boyd2004convex} holds, hence strong duality.
\end{proof}

\begin{proposition}\label{prop:equality}
    Let $v_\infty^\star, v'_\infty$ be the optimal values of the model-based counterparts of~\eqref{eq:LP-dual} and~\eqref{eq:LP-primal} respectively (corresponding to $\cG(x) = \{G(x)\}$), and assume that~\eqref{eq:OCP} has a finite value. If for $x\in\cX$, $\ell(x,\cdot)$ is convex over $\cU$, then $v_\infty^\star = v'_\infty = V^\star$.
\end{proposition}
\begin{proof}
    This follows from either of~\cite[Theorems 2.1 \& 2.3]{vinter1993convex} and~\cite[Theorem 2.3]{lasserre2008nonlinear}. Defining $$F(x) \triangleq \left\{ f(x) + G(x) \, u \mid u\in\cU \right\}, \qquad \ell^\star(x,y) \triangleq \min_{u\in\cU} \left\{\ell(x,u)\mid y = G(x) \, u \right\},$$
    we verify the hypotheses of the invoked theorems:
    \begin{itemize}
        \item $\ell^\star$ and $h$ are lower semi-continuous, because $\ell,h,f$ and $G$ are continuous on their respective sets and $\cU$ is compact;
        \item $[0,T]$ and $\cX$ are compact;
        \item For all $x\in\cX$, $F(x)$ is an affine image of $\cU$ (compact, convex set), hence it is compact and convex, and the set $\left\{ (t,x,y) \mid t\in [0,T], \ x\in\cX, \ f(x) + y \in F(x) \right\}$ is compact;
        \item \eqref{eq:OCP} is assumed to have a finite solution;
        \item $\ell^\star(x,\cdot)$ is convex on $\R^n$, because it is a minimization of a convex function ($\ell(x,\cdot)$ is convex by assumption) with a linear constraint, which preserves convexity;
        \item Strong duality $v^\star_\infty = v'_\infty$ holds by Proposition~\ref{prop:strong_duality}.
    \end{itemize}
    Since the set of polynomials $\R[t,x]$ ($(t,x)\in[0,T]\times\cX$) is dense in $C^1\left([0,T]\times\cX\right)$ with respect to the $C^1$ norm~\cite{nachbin1949algebres}, we conclude the proof.
\end{proof}

\begin{remark}
    Proposition~\ref{prop:equality} cannot be extended to the data-driven subsolution $w_\infty^\star$. This is due to the bilinearity of $(B,u)\mapsto B\,u$, which does not guarantee the convexity of $\tilde{F}(x) \triangleq \left\{f(x) + B \, u \mid (B,u) \in\cG(x)\times\cU \right\}$ anymore.
\end{remark}

\subsection{Convergence with exciting datasets}

\begin{definition}\label{def:data_excitation}
    For $\epsilon,c>0$, the dataset $\cD$ is said to be $(\epsilon,c)$-exciting if, for all $x\in\cX$, there exist $m$ indices $\{i_1,\cdots,i_m\}$ of data points such that $\|x-x_{i_j}\|\leq\epsilon$ for all $j\in[m]$, and $\rho_{\min}(\hat{U})\geq c$, with $\hat{U} \triangleq \left[\frac{u_{i_1}}{\|u_{i_1}\|},\cdots,\frac{u_{i_m}}{\|u_{i_m}\|}\right]$ and $\rho_{\min}(\cdot)$ giving the minimal singular value of a matrix.
\end{definition}

\begin{lemma}
    If $\cD$ is $(\epsilon,c)$-exciting, then
    \begin{equation}
        \|B-G(x)\|\leq\epsilon' \triangleq \frac{2\epsilon K\sqrt{m}}{c} \qquad \forall(x,B)\in\Gamma. \label{eq:lemma}
    \end{equation}
\end{lemma}
\begin{proof}
    For a $j\in[m]$ and $x\in\cX$, we have
    \begin{align*}
        \|(B-G(x))u_{i_j}\| \leq \|Bu_{i_j}-y_{i_j}\|+\|(G(x_{i_j})-G(x))u_{i_j}\| \leq 2K\|u_{i_j}\|\|x-x_{i_j}\| \leq 2K\|u_{i_j}\|\epsilon
    \end{align*}
    and $\|(B-G(x))\hat{U}\| \leq \|(B-G(x))\hat{U}\|_F$, with
    \begin{align*}
        \|(B-G(x))\hat{U}\|_F &= \sqrt{\sum_{j=1}^m \left\|(B-G(x))\frac{u_{i_j}}{\|u_{i_j}\|}\right\|^2} \leq \sqrt{m\max_{j\in[m]}\left\|(B-G(x))\frac{u_{i_j}}{\|u_{i_j}\|}\right\|^2} \\
        &\leq \sqrt{m}\max_{j\in[m]}\frac{\|(B-G(x))u_{i_j}\|}{\|u_{i_j}\|}.
    \end{align*}
    Hence $\|(B-G(x))\hat{U}\| \leq 2 \epsilon K\sqrt{m}$. Since $\hat{U}$ is square ($m$ columns of vectors in $\R^m$) and invertible ($\rho_{\min}(\!\hat{U}\!)\!>\!0$), then 
    \begin{align*}
        \|(B-G(x))\hat{U}\,\hat{U}^{-1}\| \leq \frac{\|(B-G(x))\hat{U}\|}{\rho_{\min}(\hat{U})} \quad \Rightarrow \quad \|B-G(x)\| \leq \frac{2\epsilon K\sqrt{m}}{c} \tag*{\qedhere}
\end{align*}
\end{proof}

\begin{theorem} \label{thm:convergence}
    For every $(\epsilon,c)$-exciting dataset $\cD$, if $\ell(x,\cdot)$ is convex over $\cU$ for all $x\in\cX$ and~\eqref{eq:OCP} has a finite optimal value, then the following is true for any feasible solution $v$ to the model-based instance of~\eqref{eq:LP-dual}
    \begin{enumerate}[i.]
        \item $v(0,x_0) - w_\infty^\star \leq \epsilon'T \,\Phi(v)$
        \item $0 \leq V^\star - w_\infty^\star \underset{\epsilon\rightarrow0}{\longrightarrow} 0$
    \end{enumerate}
    with $\Phi: v\mapsto\underset{\cU\times[0,T]\times\cX}{\max} \|u\| \|\frac{\partial v}{\partial x}(t,x)\|$ and $\epsilon'$ as in~\eqref{eq:lemma}.
\end{theorem}
\begin{proof}
    For $v$ a feasible solution to the model-based instance of problem~\eqref{eq:LP-dual}, \emph{i.e.} $\cG(x) = \{G(x)\}$, the modified candidate $\hat{w}:(t,x) \mapsto v(t,x)+\epsilon'\Phi(v)(t-T)$ is feasible for the data-driven instance of~\eqref{eq:LP-dual}. Indeed, let $(t,x,B,u)\in[0,T]\times\Gamma\times\cU$:
    \begin{align*}
        \begin{split}
            \frac{\partial \hat{w}}{\partial t} + \frac{\partial \hat{w}}{\partial x} \Big(f(x) + B \, u\Big) + \ell(x,u) &= \frac{\partial v}{\partial t} + \epsilon'\Phi(v) + \frac{\partial v}{\partial x} \Big(f(x) + B \, u\Big) + \ell(x,u) \\
            &= \frac{\partial v}{\partial t} + \frac{\partial v}{\partial x} \Big( f(x) + G(x) \, u \Big) + \ell(x,u)  + \frac{\partial v}{\partial x} \Big( B - G(x) \Big)u + \epsilon'\Phi(v) \\
            &\overset{\eqref{eq:subHJB}}{\geq} \frac{\partial v}{\partial x} \Big( B - G(x) \Big)u + \epsilon'\Phi(v) \\
            &\geq - \left\|\frac{\partial v}{\partial x} \Big( B - G(x) \Big)u \right\| + \epsilon'\Phi(v) \\
            &\overset{\eqref{eq:lemma}}{\geq} - \epsilon'\Phi(v) + \epsilon'\Phi(v) = 0
        \end{split}
    \end{align*}
    and $\hat{w}(T,x) = v(T,x) \overset{\eqref{eq:subTerminal}}{\leq} h(x)$. This gives $\hat{w}(0,x_0) = v(0,x_0)-\epsilon'T\,\Phi(v) \leq \sup_w w(0,x_0) = w_\infty^\star$.
    As for (ii), we have $w_\infty^\star\leq V^\star$ from Proposition~\ref{prop:inequalities}. Moreover, $V^\star = \sup_v v(0,x_0)$ according to Proposition~\ref{prop:equality} under the theorem's assumptions. So, for any arbitrary $\eta>0$, there exists a feasible $v_\eta$ such that $V^\star-\eta\leq v_\eta(0,x_0)$. This gives $V^\star - w^\star_\infty \leq \eta + v_\eta(0,x_0) - w^\star_\infty \leq \eta + \epsilon' T \Phi(v_\eta) \underset{\epsilon'\to 0}{\longrightarrow} \eta$. Since $\eta > 0$ was arbitrary, the sandwich rule concludes.
\end{proof}

\section{Moment-SoS hierarchy framework} \label{sec:msos}

The moment-\ac{SoS} hierarchy leverages Putinar's Positivstellensatz~\cite{Putinar93} to reformulate sign constraints as follows: for $\xi = (t,x,B,u) \in \Xi$, defining the \ac{HJB} slack
$$ S(\xi) = \frac{\partial w}{\partial t}(t,x) \! + \! \frac{\partial w}{\partial x}(t,x) \Big(f(x) \! + \! B \, u\Big) + \ell(x,u), $$
\eqref{eq:subHJB} is strengthened as $S \in \cQ(\Xi)$ and \eqref{eq:subTerminal} is recast as $h - w(T,\cdot) \in \cQ(\cX_T)$, with quadratic modules defined as
\begin{equation} \label{eq:QmodTerm}
    \cQ(\cX_T) \triangleq \left\{\sigma + \sigma_T^\top p_T \mid \sigma, \sigma_T \ \text{sums of squares}\right\}
\end{equation}
$$\cQ(\Xi) = \left\{\sigma + \underset{\triangleq \, p_0(t)}{\underbrace{(T-t) t}} \; \sigma_0 + \sigma_\cX^\top p_\cX + \sum_{i=1}^N \sigma_i p_i + \sigma_\cU^\top p_\cU \right\}.$$
Then, the moment-\ac{SoS} hierarchy consists in restricting the search for the sums of squares $\sigma_\bullet$ to polynomials of degree at most $d \in \N$, which is called \ac{SoS} strengthening. By Lagrange duality, on the dual side, this amounts to replacing the measures $\mu$ and $\nu$ by vectors stacking their moments up to degree $d$, which we denote by $z$ and $\vartheta$, respectively. The dual to degree $d$ quadratic module membership then requires the Riesz functional  {introduced next.}
\begin{definition} \label{def:riesz}
    Let $\chi_d$ (resp. $\psi_d$) be the basis of $\R_d[x]$ (resp. $\R_d[\xi]$) in which we compute moments, and $Q(x) = q^\top \chi_d$ (resp. $S(\xi) = s^\top\psi_d$) be a polynomial in $x \in \cX$ (resp. $\xi \in \Xi$). The Riesz functional associated with $\vartheta$ (resp. $z$) is:
    $$ L_\vartheta(Q) = q^\top \vartheta , \quad \text{resp.} \quad L_z(S) = s^\top z. $$
    In particular, extending this definition to vectors and matrices of polynomials, $\vartheta = L_\vartheta(\chi)$ (resp. $z = L_z(\psi)$).
\end{definition}
With this definition, introducing $k=\lfloor \nicefrac{d}{2}\rfloor$ for notational brevity, the hierarchy provides~\eqref{eq:LP-primal} with moment relaxations:
\begin{subequations}\label{eq:moments}
\begin{gather}
    \hspace{-2em} w_d' \triangleq \inf_{z,\vartheta} L_z(\ell) + L_\vartheta(h) \\
    \text{s.t. } \ \forall w \in \R_{d-\deg(f+Bu)+1}[t,x], \nonumber \\
    \hspace{-1em} L_\vartheta(w(T,\cdot)) - w(0,x_0) = L_z\left(\frac{\partial w}{\partial t} + (f+B \, u) \frac{\partial w}{\partial x}\right) \label{eq:momliouville} \\
    L_\vartheta\Big(\chi_k \chi_k^\top \Big) \succeq 0, \quad L_\vartheta \Big(p_T \chi_{k_T}\chi_{k_T}^\top \Big) \succeq 0 \label{eq:momterm} \\
    L_z\Big(\psi_k\psi_k^\top \Big) \succeq 0, \quad L_z \Big(p_\bullet \psi_{k_\bullet} \psi_{k_\bullet}^\top \Big) \succeq 0 \label{eq:momoccup} \\
    \bullet \in \{0,\cX,\cU,1,\ldots,N\}, \quad k_{T/\bullet} = \lfloor \nicefrac{d - \deg p_{T/\bullet}}{2} \rfloor \nonumber
\end{gather}
\end{subequations}
Importantly, $L_z(\psi_k\psi_k^\top)$ is called the order $d$ \textit{moment matrix} of $z$, and $L_\vartheta(p_T \chi_{k_T}\chi_{k_T}^\top)$ is called the order $d$ \textit{localizing matrix} of $\vartheta$. The \ac{SoS} multipliers $\sigma$ and $\sigma_T$ (resp. $\sigma,\sigma_\bullet$) in the definition~\eqref{eq:QmodTerm} of $\cQ(\cX_T)$ (resp. $\cQ(\Xi)$) are the Lagrange multipliers of \ac{LMI} constraints~\eqref{eq:momterm} (resp.~\eqref{eq:momoccup}). 

\begin{proposition} \label{prop:momconv}
    Let $w^\star_d$ denote the optimal value of the \ac{SoS}-strengthening of~\eqref{eq:LP-dual}. Then, the hierarchy converges monotonically, \emph{i.e.}
    \begin{equation} \label{eq:convergence}
        w'_d \underset{d\to\infty}{\nearrow} w^\star_\infty.
    \end{equation}
    and strong duality $w'_d = w^\star_d$ holds at each relaxation order.
\end{proposition}
\begin{proof}
    The proof of monotone convergence is identical to~\cite[Theorem 3.6]{lasserre2008nonlinear}. By~\cite[Proposition 6]{Tacchi2022convergence}, it only remains to prove that moment relaxations~\eqref{eq:moments} of~\eqref{eq:LP-primal} enforce uniform boundedness of the pseudo-moments of order $0$. Let $z,\vartheta$ be pseudo-moment sequences corresponding to $\mu,\nu$ respectively. Then, testing~\eqref{eq:momliouville} against $w(t,x)=1$ yields
    \begin{align*}
        \vartheta_0 & =L_\vartheta(1) = L_\vartheta\Big( w(T,\cdot)\Big) \\
        & = w(0,x_0) + L_z\left( \cancel{\frac{\partial w}{\partial t}} + (f+B \, u) \cancel{\frac{\partial w}{\partial x}} \right) \\
        & = w(0,x_0) = 1
    \end{align*}
    and testing it against $w(t,x) = t$ yields
    \begin{align*}
        z_0 & = L_z(1) = L_z\left(\frac{\partial w}{\partial t}\right) = L_z\left( \frac{\partial w}{\partial t} +  \cancel{\frac{\partial w}{\partial x}} (f+B \, u) \right) \\
        & = L_\vartheta\big(w(T,\cdot)\big) - \cancel{w(0,x_0)} = T.
    \end{align*}
    As a result, when $d$ is large enough, $\vartheta_0=1$ and $z_0=T$ are fixed, hence uniformly bounded, and applying~\cite[Proposition 6]{Tacchi2022convergence} concludes the proof.
\end{proof}

\section{Extracting a control law} \label{sec:extract}

Throughout, $w_d$ denotes an optimal solution of the degree $d$ \ac{SoS} strengthening of~\eqref{eq:LP-dual} from section~\ref{sec:msos}, $\mu$ the optimal occupation measure of~\eqref{eq:LP-primal} on $\Xi = [0,T]\times \Gamma \times \cU$, in the variables $\xi = \left(t,(x,B),u\right)$, and $z_d$ the corresponding optimal pseudo-moment sequence in the degree $d$ moment relaxation.

\subsection{From the dual on functions}

A standard way of extracting a control law was proposed by Henrion et al. in~\cite{henrion2008nonlinear} (in their case, the LP was model-based, i.e. $\cG(x) = \{G(x)\}$): one can define the \ac{HJB} slack, for $\xi = \left(t,(x,B),u\right) \in \Xi$, as
$$ S_d(\xi) = \frac{\partial w_d}{\partial t}(t,x) + \frac{\partial w_d}{\partial x}(t,x)\Big(f(x) + B \, u \Big) + \ell(x,u), $$
which by design is a nonnegative polynomial on $\Xi$. In the infinite-dimensional setting, if an optimal $w_\infty$ existed, complementary slackness would enforce
$ \int S_\infty \;  d\mu = 0 $, hence $S_\infty$ should vanish on the support of $\mu$, which encodes the trajectories. This justifies the heuristic computation of a feedback control, at given $(t,x) \in [0,T]\times \cX$, in the form
\begin{equation} \label{eq:feedbackHJB}
    \kappa_d(t,x) = \argmin{u \in \cU} \min_{B \in \cG(x)} S_d(t,x,B,u).
\end{equation}
Since $S_d$ is affine in $B$ and the $\cG(x)$ are intersections of convex quadratic sets, the inner minimization is a small convex program. If $\ell(x,u)$ is not strictly convex in $u$, the argmin may be set-valued, in which case any measurable selection can be used. The feedback~\eqref{eq:feedbackHJB} follows the optimistic plan exactly: it would be optimal in the best-case scenario where $G(x)$ is the most favorable. On the true system, this is a source of bias, which we avoid using the next extraction method.

\subsection{From the primal on moments}

Here we propose an alternative, averaging a $(t,x)$ disintegration of the occupation measure, over all possible models $B$ and controls $u$:
$$ \bar\kappa(t,x) = \E_\mu[u \mid t,x] = \int u \; d\mu(u,B \mid t,x), $$
defined $\mu^{(t,x)}$-a.e. and with values in $\cU$ since $\cU$ is convex. By design, this feedback policy averages all inputs and models charged by $\mu$, so that it does not only account for the best-case scenario (although $\mu$ is optimized for such a scenario). By definition of the conditional mean, one can view $\bar\kappa$ as the orthogonal projection, in $L^2(\mu;\R^m)$, of the projection function $\pi_u : (t,x,B,u) \longmapsto u$, onto the subspace $L^2(\mu^{(t,x)};\R^m)$ of functions depending only on $(t,x)$. Its best polynomial approximation of degree at most $r$ in $L^2(\mu^{(t,x)};\R^m)$ directly is the orthogonal projection of $\pi_u$ onto $\R_r[t,x]$, since by Pythagoras it holds for any such polynomial $\kappa_r$
$$ \|\kappa_r - \pi_u\|^2_{L^2(\mu)} = \underset{\in L^2(\mu^{(t,x)};\R^m)}{\|\underbrace{\kappa_r - \bar\kappa}\|^2_{L^2(\mu)}} + \underset{\in L^2(\mu^{(t,x)};\R^m)^\perp}{\|\underbrace{\bar{\kappa} - \pi_u}\|^2_{L^2(\mu)}}. $$
Since in practice only pseudo-moment vectors $\vartheta$ and $z$ up to degree $d \geq r$ are available, denoting $\varphi_r(t,x)$ a size $\omega_r = \binom{r+n+1}{r}$ basis of $\R_r[t,x]$, we look for a matrix $A \in \R^{\omega_r \times m}$ to parameterize the representation
$$ \kappa_r(t,x) = A^\top \varphi_r(t,x), $$
and the $L^2(\mu)$-norm to be minimized is recast as
\begin{align*}
    L_z(\|\kappa_r - \pi_u\|^2) & = L_z\left(\|A^\top \varphi_r(t,x) - u\|^2\right) \\
    & \hspace{-6.91pt} = \tr\left(A^\top M_r \, A\right) - 2 \tr \left(A^\top C_r \right) + L_z(\|u\|^2), 
\end{align*}
$M_r = L_z(\varphi_r\varphi_r^\top) \in \R^{\omega_r \times \omega_r}$, $C_r = L_z(\varphi_r \, u^\top) \in \R^{\omega_r \times m}$. 

The corresponding optimality condition
is
\begin{equation} \label{eq:KKT-cond}
    M_r \, \hat{A} = C_r,
\end{equation}
which recovers the moment-matching extraction of~\cite{majumdar2014convex}, where a polynomial controller is obtained from the moments of a control measure that corresponds to $(u \, \mu)^{(t,x)} = \bar{\kappa} \, \mu^{(t,x)}$ in our framework. The novelty here is that our projection viewpoint allows us to add input constraints.

Indeed, while $\bar\kappa$ takes values in $\cU$, such a guarantee disappears when we use polynomial approximation, and overshoots occur when the input saturates over a region, even on the support of the occupation measure. Therefore, we restrict the approximation to polynomials with values in $\cU$, by solving the problem
\begin{equation} \label{eq:extract}
\begin{aligned}
\hat{A} \in \argmin{A \in \R^{\omega_r \times m}} \tr\left(A^\top M_r^z \, A \right) - 2 \, \tr\left(A^\top C_r^z\right) \\
\text{s.t.} \quad \forall (t,x) \in [0,T]\times \cX, \quad A^\top \varphi_r(t,x) \in \cU \\
\textit{i.e.} \quad P_\cU \, A^\top \varphi_r(t,x) - q_\cU \geq 0.
\end{aligned}
\end{equation}
Since $\cU = \{u \in \R^m \mid P_\cU \, u - q_\cU \geq 0\}$ is a polytope, the constraint is pointwise linear in $A$ and can be represented through Putinar certificates with fixed degree. $M_r^z$ being positive semidefinite by design,~\eqref{eq:extract} can then be recast as a small semidefinite program. We call the obtained $\hat\kappa_r = \hat{A}^\top \varphi_r(t,x)$ the \ac{CD} feedback, because it relies on an $L^2(\mu)$ orthogonal projection approach.

\begin{remark}
    With an initialization in $x_0$, i.e. with $\delta_{(0,x_0)}$ in the LHS of~\eqref{eq:liouville}, the optimal $\mu^{(t,x)}$ from~\eqref{eq:LP-primal} is supported on a curve, so that $\bar\kappa$ is only defined along that curve, and the computed $M_r$ is ill-conditioned because it aims at a singular moment matrix. For the extraction, we therefore solve the relaxation with the initial condition chosen uniformly at random on a small box around $x_0$ (in our experiments $x_0 + [-0.1,0.1]^n$, which can be interpreted as 
    a bounded 
    uncertainty on the initial condition). Then, the optimal $\mu^{(t,x)}$ is supported on a tube on which $\bar\kappa$ is defined, and the optimal moment matrix is positive definite. Outside the tube, only the polynomial approximation $\hat\kappa_r = \hat{A}^\top \varphi_r(t,x)$ is defined as an extrapolation. The bound $w_d^\star$ remains computed with initialization at $x_0$; with bounded additive noise $\Delta x$, the optimal value bounds $\E_{\Delta x}[V(0,x_0 + \Delta x)]$ instead.
\end{remark}

\section{Numerical application}\label{sec:num_app}

\subsection{Choosing data location}

While the number of \ac{LMI} in~\ref{eq:moments} grows linearly with the number $N$ of data points, the computational time still grows sharply with $N$ (because the linearity coefficient itself increases combinatorially in the relaxation degree $d$). As a result, one might want to minimize the number of data points by using only well-localized samples, either from existing data or, when possible, from new experiments. For a computed optimal solution $w_d(t,x)$ to the degree $d$ \ac{SoS} strengthening of~\eqref{eq:LP-dual}, we define the state score
$$ S_d^\star(x) \triangleq \min_{\substack{t \in [0,T], u \in \cU \\
B \in \cG(x)}} S_d(t,x,B,u) \geq 0.$$
By definition, $S_d^\star$ is small where the \ac{HJB} inequality is close to saturation, which corresponds to the support of the optimal occupation measure $\mu$, so that small values of $S_d^\star$ localize states where the bound is decided, and drawing data there tightens $\cG(x)$ where it is most likely to improve. For instance, new samples can be drawn without replacement among some candidates with probabilities proportional to $\exp(-\nicefrac{S_d^\star(x)}{\tau})$, where $\tau = 0.1\,\mathrm{median}(S^\star_d)$, and the median is empirically computed over a finite set of candidates.

From this observation, two approaches are possible, depending on whether one has access only to a large pool of measurements or to a sampling experiment. This paper focuses on the case where no simulator or experimental benchmark is available; one can rely on a fixed pool of $P = 2000$ samples (uniform on $\cX \times \cU)$. Since the cost of the update grows with $N$, accounting for all the data might be too expensive, so one can start with a small number of samples (10 in our experiments), and add 10 per round, with various policies for the selection of candidates:
\begin{itemize}
    \item one can pick the new data points uniformly at random among the available measurements;
    \item one can choose them by greedy coverage, \emph{i.e.} taking the farthest point in $x$ from the current dataset;
    \item one can leverage the slack score $S^\star_d$, as mentioned above, to select new data points.
\end{itemize}
Table~\ref{tab:selection} illustrates how the data selection process affects the methodology's overall accuracy using a simple example: the Van der Pol oscillator (see Section~\ref{sec:benchmarks}). Each round solves one degree $d=4$ relaxation~\eqref{eq:moments} (in about 6s), and the results show that slack-based selection with 40 samples gives a bound comparable to the one obtained with 200 random samples, as well as a worst-case cost (4.17) much better than with 40 random samples (4.51) and close to that with 200 (4.09). Coverage improves less than slack: what matters is information close to the optimal trajectories, not filling $\cX$, which is only instrumental in the proof of Theorem~\ref{thm:convergence}.

\begin{table}[ht]
\centering
\caption{\footnotesize Selection from a pool of 2000 samples (Van der Pol, $d=4$, mean $\pm$ standard
error over 10 seeds; $V^\star=3.746$). True: cost of the \ac{CD} feedback on the true system; worst: over the
consistent gains of section~\ref{sec:benchmarks}.}\label{tab:selection}
\begin{tabular}{@{}llccc@{}}
\hline
$N$ & rule & $w^\star_d$ & true & worst\\
\hline
20 & random   & $2.857\pm0.054$ & $3.900\pm0.009$ & $4.92\pm0.21$\\
20 & coverage & $2.927\pm0.046$ & $3.908\pm0.020$ & $4.65\pm0.08$\\
20 & slack    & $3.067\pm0.034$ & $3.840\pm0.004$ & $4.28\pm0.04$\\
40 & random   & $2.987\pm0.037$ & $3.878\pm0.012$ & $4.51\pm0.08$\\
40 & coverage & $3.041\pm0.020$ & $3.852\pm0.007$ & $4.36\pm0.06$\\
40 & slack    & $3.158\pm0.008$ & $3.846\pm0.004$ & $4.17\pm0.04$\\
200 & random  & $3.152\pm0.010$ & $3.842\pm0.004$ & $4.09\pm0.03$\\
\hline
\end{tabular}
\end{table}

If new data points can be obtained at will, then it is possible to investigate active sampling, \emph{i.e.} choosing the exact location of each new data point, either uniformly at random, or using the slack score again, or even on the trajectory induced by the current feedback on the true system. We leave such sampling methods and their analysis for future work.

\subsection{Benchmarks and methodology} \label{sec:benchmarks}

In all benchmarks we assume $G(x) = (0,\ldots,0,g(x))^\top$ with unknown gain $g$ and hence represent $B = (0,\ldots,0,b)$ with $b \in \cG(x) \subset \R$. We also work with $\cX_T = \cX$.

We first study 3 variants of a Duffing oscillator
$$ \dot{x}_1 = x_2, \quad \dot{x_2} = x_3 -x_1-x_1^3-0.3x_2, \quad \dot{x}_3 = g(x) \, u - x_3 $$
with $\cX = [-1.5,1.5]^3$, $\cU = [-1,1]$, $x_0 = (1,0,0)$, $T=6$, 
$$\ell(x,u) = x_1^2+x_2^2+0.1x_3^2+0.5u^2, \quad h(x) = x_1^2+x_2^2, \quad \text{and:}$$
\begin{itemize}[leftmargin=*]
    \item polynomial gain $g(x) = 1.2 - 0.3x_1^2 + 0.15x_2$ ($K=0.92$)
    \item non-polynomial $g(x) = 1+0.4\sin x_1 \cos x_2$ ($K=0.57$)
    \item noisy $\dot{x}_3$ measurements ($\varepsilon_{\dot{x}}=0.1$) on the polynomial gain.
\end{itemize}
The second benchmark is the Van der Pol oscillator
$$ \dot{x}_1 = x_2, \quad \dot{x}_2 = (1-x_1)^2x_2 - x_1 + (1+0.4x_1^2)u \quad (K=1.6) $$
with $\cX = [-2,2]^2$, $\cU = [-1,1]$, $x_0 = (1,1)$, $h(x) = \|x\|^2$, $\ell(x,u) = \|x\|^2 + u^2$ and $T=2$. 
The reference optimal values $V^\star$ (Duffing: 1.9988 for polynomial gain, 1.7958 for non-polynomial; Van der Pol: 3.7456) 
are computed by direct shooting with the true gain (polynomial gains are confirmed to be tight using the methodology of~\cite{lasserre2008nonlinear}).

We solve the moment-\ac{SoS} relaxations of degree $d=4$, using time splitting~\cite{CKH22} into three pieces of $2$s (Duffing, $T=6$), in order to avoid bad \ac{LMI} conditioning. The polynomial \ac{CD} feedback $\hat{\kappa}_r$ is fitted on each piece and has degree $r=2$. 
All SDPs are solved with MOSEK via YALMIP in MATLAB on an Apple M2 laptop (16GB).

We compare our method with two existing MPC schemes that use the same information to identify a model of $G$: LI-MPC uses Lipschitz interpolation, \emph{i.e.}, the midpoint of $\cG(x)$ (kinky inference~\cite{Cal17}, nominal version of~\cite{MLMC20}); SINDy-MPC uses SINDy~\cite{BPK16,KKB18},
a sparse regression of $g$ on the monomials of degree at most 3. Once these model candidates are learned, a standard MPC scheme is applied. We implement the MPC input as a piecewise-constant policy on a grid with step sizes of 0.2s (Duffing) or 0.1s (Van der Pol), re-planning every two steps. We generate a pool of 2000 samples uniformly at random using 10 fixed seeds, then apply our data selection method. The baselines have no access to $S^\star_d$, so they add random samples, which we also use in our method. We report the gap to $V^\star$ on the true system, and the worst cost over five gains $(1-\theta)\underline{g}+\theta\overline{g}$, where $\cG(x) = [\underline{g}(x),\overline{g}(x)]$, $\theta \in \{0, \frac{1}{4}, \frac{1}{2}, \frac{3}{4},  1\}$, a lower estimate of the worst case.

\subsection{Results and comparisons} \label{sec:results}

Table~\ref{tab:main} compares our method with the $\hat{\kappa}_2$ control policy to the baselines LI-MPC and SINDy-MPC. For all tested benchmarks, MPC performs best against the true system, while our method appears more robust in the worst-case scenario when data are selected according to the slack score, which the baselines cannot use. Sample selection also reduces the $V^\star - w_d^\star$ gap.

\begin{table}[t]
\caption{\footnotesize Closed-loop performance at $N=40$, mean $\pm$ standard error over 10 seeds: gap to $V^\star$ on the true system, worst cost over consistent gains. Best per benchmark in bold.} \label{tab:main}
\centering\scriptsize\setlength{\tabcolsep}{3pt}
\begin{tabular}{@{}llcc@{}}
\hline
benchmark & controller & true (\%) & worst \\
\hline
Duffing, base & ours, selected & $+2.3\pm.1$ & $\mathbf{2.31\pm.02}$ \\
 & ours, random & $+3.5\pm.4$ & $2.80\pm.05$ \\
 & LI-MPC & $+0.2\pm.0$ & $2.67\pm.05$ \\
 & SINDy-MPC & $\mathbf{+0.1\pm.0}$ & $2.67\pm.05$ \\
\hline
Duffing, non-poly. & ours, selected & $+1.3\pm.1$ & $\mathbf{1.96\pm.01}$ \\
 & ours, random & $+2.2\pm.1$ & $2.14\pm.02$ \\
 & LI-MPC & $\mathbf{+0.1\pm.0}$ & $2.02\pm.02$ \\
 & SINDy-MPC & $\mathbf{+0.1\pm.0}$ & $2.02\pm.02$ \\
\hline
Duffing, noise & ours, selected & $+4.7\pm.3$ & $\mathbf{2.62\pm.02}$ \\
 & ours, random & $+7.1\pm.7$ & $3.26\pm.08$ \\
 & LI-MPC & $+0.2\pm.0$ & $3.20\pm.12$ \\
 & SINDy-MPC & $\mathbf{+0.1\pm.0}$ & $3.21\pm.12$ \\
\hline
Van der Pol, $T=2$ & ours, selected & $+2.7\pm.1$ & $\mathbf{4.17\pm.04}$ \\
 & ours, random & $+3.5\pm.3$ & $4.51\pm.08$ \\
 & LI-MPC & $+0.3\pm.1$ & $4.76\pm.17$ \\
 & SINDy-MPC & $\mathbf{+0.1\pm.0}$ & $4.69\pm.16$ \\
\hline
\end{tabular}
\end{table}

\begin{figure}[t]
\centering
\includegraphics[width=.65\columnwidth]{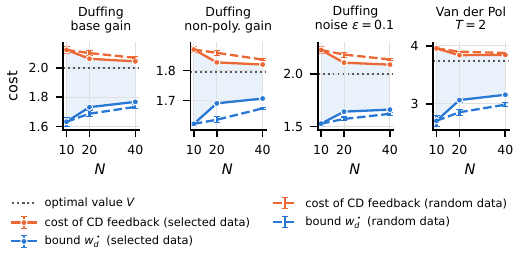}
\caption{\footnotesize Bound $w^\star_d$ (blue) and true cost of the \ac{CD} feedback (orange) against $N$, mean $\pm$ standard error over 10 seeds. Solid: selected data. Dashed: random data, as for the
baselines. Dotted: $V=V^\star$, which lies in the shaded band between the bound and the cost of our feedback.}\label{fig:sandwich}
\end{figure}

Figure~\ref{fig:sandwich} illustrates the inequality $w^\star_d \leq V^\star \leq J(\kappa)$ of Propositions~\ref{prop:inequalities} and~\ref{prop:momconv}, regardless of the choice for the feedback law $u(t) = \kappa(t,x(t))$, the MPC baselines included. Selection raises the bound at every $N$. At $N=40$, for our feedback, the certified gap $\nicefrac{J(\kappa)-w^\star_d}{J(\kappa)}$, an upper bound on $\nicefrac{J(\kappa)-V^\star}{J(\kappa)}$, is then on average 13.5\% (base Duffing), 6.2\% (non-polynomial), 20.5\% (noisy) and 17.9\% (Van der Pol), while for LI-MPC it is 11.7\%, 5.1\%, 17.0\% and 15.9\%, respectively. It is limited by the looseness of the bound (Theorem~\ref{thm:convergence}). Neither baseline provides any such certificate.

\begin{figure}[t]
\centering
\includegraphics[scale=1]{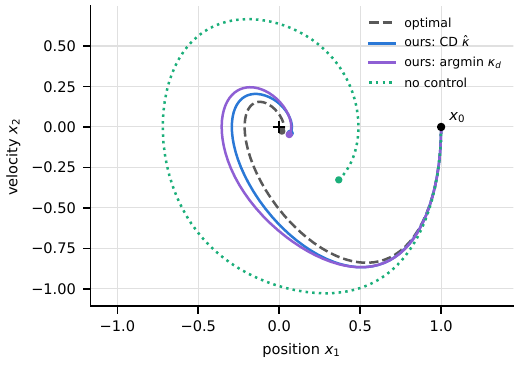}
\caption{\footnotesize Closed-loop trajectories from $x_0=(1,0,0)$ on the true system, projected on
$(x_1,x_2)$: Duffing, base gain, $N=40$ selected samples, with the optimal control (true
gain known) and without control. Dots: $x(T)$. Costs: $2.039$ ($\hat\kappa_r$), $2.107$
($\kappa_d$), $1.999$ (optimal), $3.898$ (no control).}\label{fig:traj}
\end{figure}

Figure~\ref{fig:traj} displays the trajectories corresponding to the feedbacks $\kappa_d(t,x) = \mathrm{argmin}_{u} \min_B \ S_d(t,x,B,u)$ and $\hat{\kappa}_r(t,x)$ obtained by \ac{CD} projection: both damp the oscillation as the optimal control does and end near the origin. They overshoot slightly more than the optimal trajectory, as they assume a favorable gain and actuate less (the actuator state $x_3$ peaks at $0.39$--$0.45$ versus $0.52$ at the optimum).

\begin{figure}[t]
\centering
\includegraphics[scale=1.1]{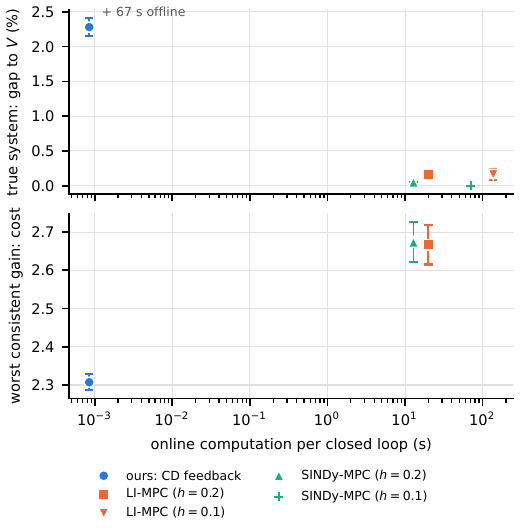}
\caption{\footnotesize Accuracy against online computation per closed loop, Duffing, base gain, $N=40$. Top: gap to true $V^\star$. Bottom: worst cost over the consistent gains. $h$: MPC step.}\label{fig:online}
\end{figure}

Our feedback is an explicit polynomial: a whole closed loop takes about 1ms in MATLAB. In contrast, the better-performing MPC schemes re-solve a nonlinear program at every re-planning step, which takes 13--156s per closed loop, \emph{i.e.} $10^4$ to $2\cdot10^5$ times more; see Figure~\ref{fig:online}. Dedicated NMPC solvers would reduce this ratio by orders of magnitude, but evaluating a polynomial will always be much faster.

\section{Conclusion}

This paper studies the optimal control of polynomial input-affine systems with an unknown but Lipschitz continuous input matrix. An uncertainty set characterizing this matrix is generated based on assumptions on the availability of a dataset and an upper bound on the Lipschitz constant. The \ac{HJB} equation used to solve the optimal control problem is then adapted to a best-case data-driven formulation of uncertainty, providing a lower bound on the best achievable cost. This equation is further relaxed into an infinite-dimensional \ac{LP} and its primal on moments.

The theoretical contributions of this work consist of a certified lower bound on the attainable cost and formal mathematical guarantees, with a convergence rate depending on the excitation properties of the data. The convergence constant grows linearly with the time horizon, the Lipschitz upper bound, and the square root of the number of inputs, illustrating the deterioration caused by large time horizons, poor Lipschitz constant estimation, or high-dimensional control inputs. A control law is extracted from the optimized primal on moments and is compared with two \ac{MPC} variants in the numerical section. The results show that it has better online computation time, as the controller is constructed offline, and with selected data, it performs better against our estimated worst case. However, it still lacks precision relative to the \ac{MPC} benchmarks, suffers from the curse of dimensionality, and remains limited to modestly sized datasets.

Future work will investigate strategies to improve the performance of our controller, for instance by sharpening the uncertainty set and model learning. The curse of dimensionality will be addressed through a combination of existing scaling strategies for SoS programming, tailored to the specific data-driven problems at hand. Likewise, we will investigate ways to mitigate the computational burden associated with large datasets while maintaining convergence properties. Eventually, we will explore possibilities to improve the controller online during operation, while avoiding both the need to rerun the whole machinery and the loss of formal guarantees.

\bibliographystyle{ieeetr}
\bibliography{RefECC27}

\end{document}